\documentclass{adm-paper}
\usepackage{amssymb,amsmath,amsthm}

\begin{document}

\title{Miniversal deformations of
chains of linear mappings\thanks{This is the authors' version of a work that was published in Algebra Discrete Math.(no.1) (2005) 47-61.}}
\author{T.N. Gaiduk,
V.V. Sergeichuk, N.A. Zharko}
\shorttitle{Miniversal
deformations}
\shortauthor{T.N. Gaiduk,
V.V. Sergeichuk, N.A. Zharko}


\communicated{V.V.
Kirichenko}

\research%

\date{}

\theoremstyle{plain}
\newtheorem{theorem}{Theorem}
\newtheorem{lemma}[theorem]{Lemma}
\newtheorem{corollary}[theorem]{Corollary}
\newtheorem{definition}[theorem]{Definition}

\renewcommand{\le}{\leqslant}
\renewcommand{\ge}{\geqslant}

\maketitle

\begin{abstract}
V.I.~Arnold [Russian Math. Surveys,
26 (no. 2), 1971, pp. 29--43] gave a
miniversal deformation of matrices
of linear operators; that is, a
simple canonical form, to which not
only a given square matrix $A$, but
also the family of all matrices
close to $A$, can be reduced by
similarity transformations smoothly
depending on the entries of
matrices. We study miniversal
deformations of quiver
representations and obtain a
miniversal deformation of matrices
of chains of linear mappings
\[
V_1 \,\frac{}{\qquad}\,
V_2\,\frac{}{\qquad}\, \cdots
\,\frac{}{\qquad}\, V_t\,,
\]
where all $V_i$ are complex or real
vector spaces and each line denotes
$\longrightarrow$ or
$\longleftarrow$.
\end{abstract}

\subjclass{2001}{15A21; 16G20}

\keywords{Parametric matrices;
Quivers; Miniversal deformations}

\section*{Introduction}

All matrices $B$ that are close to a
given square complex matrix $A$
reduce by similarity transformations
to their Jordan canonical forms, but
these forms and transformations may
be discontinuous relative to the
entries of $B$. Arnold \cite{arn}
(see also \cite[\S\,30]{arn3})
constructed a normal form, to which
not only the matrix $A$, but all
matrices close to it, can be reduced
by smooth similarity
transformations. He called this
normal form a \emph{miniversal
deformation} of $A$. A miniversal
deformation of real matrices for
similarity was given by Galin
\cite{gal}. A miniversal deformation
of pairs of $m$-by-$n$ matrices with
respect to simultaneous equivalence
(that is, of matrix pencils) was
obtained in the paper \cite{kag},
which was awarded by the SIAG/LA
(SIAM Activity Group on Linear
Algebra) Prize in Applied Linear
Algebra for the years 1997--2000.
The miniversal deformations from
\cite{kag} and \cite{gal} were
simplified significantly in
\cite{gar_ser}. These results are
important for applications in which
one has matrices that arise from
physical measurements, which means
that their entries are known only
approximately.

The notion of a miniversal deformation
was extended to quiver representations
in \cite{gar_ser}. Recall that a
\emph{quiver} is a directed graph, its
\emph{representation ${\cal A}$} over a
field $\mathbb F$ is given by assigning
to each vertex $p$ a finite dimensional
vector space ${\cal A}_p$ over $\mathbb
F$ and to each arrow $\alpha\colon p
\to q$ a linear mapping ${\cal
A}_{\alpha}\colon {\cal A}_p \to {\cal
A}_q$. Studying the family of quiver
representations whose matrices are
close to the matrices of a given
representation ${\cal A}$, we can
independently reduce the matrices of
each representation to Belitskii's
canonical form \cite{ser} losing the
smoothness relative to the entries of
these matrices. This leads to the
problem of constructing a simple normal
form to which all representations close
to ${\cal A}$ can be reduced by smooth
changes of bases; that is, to the
problem of constructing a miniversal
deformation of ${\cal A}$.

In Section \ref{s1} we recall a
theorem from \cite{gar_ser} that
admits to construct miniversal
deformations of quiver
representations. In Section \ref{s3}
we give a direct and constructive
proof of this theorem (it was
deduced in \cite{gar_ser} from some
result about miniversal deformations
formulated in \cite{arn2}). In
Section \ref{s1} we also prove that
a miniversal deformation of each
quiver representation is easily
constructed from miniversal
deformations of direct sums of two
indecomposable representations. In
Section \ref{s2} we obtain a
miniversal deformation of each
quiver $1-2-\dots-t$ with an
arbitrary orientation of its arrows;
that is, a miniversal deformation of
matrices of chains of linear
mappings
\[
V_1 \,\frac{}{\qquad}\,
V_2\,\frac{}{\qquad}\, \cdots
\,\frac{}{\qquad}\, V_t\,,
\]
where all $V_i$ are complex or real
vector spaces and each line denotes
$\longrightarrow$ or
$\longleftarrow$.

\section{Miniversal deformations of
quiver representations}\label{s1}

In this section we consider
representations of a quiver $Q$ with
vertices $1,\dots,t$. Let ${\cal A}$ be
any representation of $Q$ over a field
$\mathbb F$. Choosing bases in the
spaces ${\cal A}_1,\dots,{\cal A}_t$ we
may give ${\cal A}$ by the matrices of
its linear mappings ${\cal
A}_{\alpha}\colon {\cal A}_p \to {\cal
A}_q$. This leads to the following
definitions. By a {\it matrix
representation} of dimension $\vec
n=(n_1,\dots,n_t)\in \{0,1,2,\dots\}^t$
of $Q$ over $\mathbb F$ we mean any set
$A$ of matrices $A_{\alpha}\in {\mathbb
F}^{\,n_q\times n_p}$ assigned to all
arrows $\alpha\colon p\to q$. Two
matrix representations $A$ and $B$ of
dimension $\vec n$ are
\emph{isomorphic} if there is a
sequence $S=(S_1,\dots,S_t)$ of
nonsingular $n_1\times
n_1,\dots,n_t\times n_t$ matrices such
that
\[
B_{\alpha}=
S_qA_{\alpha}S_p^{-1}\qquad
\text{for each arrow $\alpha\colon
p\to q$.}
\]
In this case we say that $S$ is an
\emph{isomorphism} of the
representations $A$ and $B$ and
write $S\colon
A\xrightarrow{\sim}B$. Clearly, all
isomorphic matrix representations
define the same (operator)
representation with respect to
different bases of its spaces.
Denote by ${\cal R}(\vec n,\mathbb
F)$ the vector space of all matrix
representations of dimension $\vec
n$ over $\mathbb F$.

From this point on, $\mathbb F$ is a
field of complex or real numbers, and
we consider only matrix representations
omitting usually the word ``matrix''
for abbreviation. A \emph{deformation}
of $A\in{\cal R}(\vec n,\mathbb F)$ is
a matrix representation  ${\cal
A}(\vec\lambda)$, $\vec\lambda
=(\lambda_1,\dots,\lambda_k)$, such
that the entries of its matrices are
convergent in a neighborhood of $\vec
0$ power series of variables (they are
called \emph{parameters})
$\lambda_1,\dots, \lambda_k$ over
$\mathbb F$ and ${\cal A}(\vec 0)=A$.
Two deformations ${\cal
A}(\vec\lambda)$ and ${\cal
B}(\vec\lambda)$ of $A\in{\cal R}(\vec
n,\mathbb F)$ are called {\it
equivalent} if the identity isomorphism
\begin{equation*}\label{bhg}
I_{\vec n}
=(I_{n_1},\dots,I_{n_t})\colon
A\xrightarrow{\sim} A
\end{equation*}
possesses a deformation ${\cal
I}(\vec\lambda)$ (its matrices are
convergent in a neighborhood of
$\vec 0$ matrix power series and
${\cal I}(\vec 0)=I_{\vec n}$) such
that
\[
{\cal B}_{\alpha}(\vec\lambda)=
{\cal I}_q(\vec\lambda) {\cal
A}_{\alpha}(\vec\lambda) {\cal
I}_p(\vec\lambda)^{-1}\qquad
\text{for each arrow $\alpha\colon
p\to q$,}
\]
in a neighborhood of $\vec 0$.

\begin{definition}\label{d}
A deformation ${\cal
A}(\lambda_1,\dots,\lambda_k)$ of a
representation $A$ is called
\emph{versal} if every deformation
${\cal B}(\mu_1,\dots,\mu_l)$ of $A$ is
equivalent to a deformation of the form
${\cal A}(\varphi_1(\vec\mu),\dots,
\varphi_k(\vec\mu)),$ where
$\varphi_i(\vec\mu)$ are convergent in
a neighborhood of $\vec 0$ power series
such that $\varphi_i(\vec 0)=0$. A
versal deformation ${\cal
A}(\lambda_1,\dots,\lambda_k)$ of $A$
is called \emph{miniversal} if there is
no versal deformation having less than
$k$ parameters.
\end{definition}

A miniversal deformation of any
representation $A$ of dimension
$\vec n$ can be constructed as
follows. The triples consisting of
all arrows $\alpha\colon p\to q$ of
$Q$ and indices of the
$n_q$-by-$n_p$ matrices $A_{\alpha}$
of $A$ form the set
\begin{equation}\label{a1}
\Upsilon_{\vec n} :=\{(\alpha,i,j)\,|\,
\alpha\colon p\to q,\ \ i=1,\dots,n_q,\
\ j=1,\dots,n_p\}.
\end{equation}
For each $(\alpha,i,j)\in\Upsilon_{\vec
n}$, define the \emph{elementary
representation} $E_{\alpha ij}$ whose
matrices are zero except for the matrix
assigned to $\alpha$; the
$(i,j)^{\text{th}}$ entry of this
matrix is $1$ and the others are $0$.

For each subset
$\varGamma\subset\Upsilon_{\vec n}$,
define the deformation
\begin{equation}\label{a2}
{\cal U}_{\varGamma}(\vec
{\varepsilon}):=
A+\sum_{(\alpha,i,j)\in\varGamma}
\varepsilon_{\alpha ij}E_{\alpha ij}
\end{equation}
of $A$, in which all
$\varepsilon_{\alpha ij}$ are
independent parameters. The
deformation
\begin{equation}\label{edr}
{\cal U}(\vec {\varepsilon}):={\cal
U}_{\Upsilon_{\vec n}}(\vec
{\varepsilon})
\end{equation}
is universal in the sense that each
deformation ${\cal
B}(\mu_1,\dots,\mu_l)$ of $A$ has
the form ${\cal U}(\vec{\varphi}
(\mu_1,\dots,\mu_l)),$ where
$\varphi_{\alpha
ij}(\mu_1,\dots,\mu_l)$ are
convergent in a neighborhood of
$\vec 0$ power series such that
$\varphi_{\alpha ij}(\vec 0)= 0$.
Hence the deformation ${\cal
B}(\mu_1,\dots,\mu_l)$ in Definition
\ref{d} can be replaced by ${\cal
U}(\vec {\varepsilon})$, which
proves the following lemma.

\begin{lemma}\label{lem}
The following two conditions are
equivalent for any deformation
${\cal
A}(\lambda_1,\dots,\lambda_k)$ of a
representation $A$:
\begin{itemize}
  \item[\rm(i)]
The deformation ${\cal
A}(\lambda_1,\dots,\lambda_k)$ is
versal.
  \item[\rm(ii)]
The deformation ${\cal U}(\vec
{\varepsilon})$ defined in
\eqref{edr} is equivalent to a
deformation of the form ${\cal
A}(\varphi_1(\vec{\varepsilon}),\dots,
\varphi_k(\vec{\varepsilon})),$
where $\varphi_i(\vec{\varepsilon})$
are convergent in a neighborhood
of\/ $\vec 0$ power series such that
$\varphi_i(\vec 0)=0$.
\end{itemize}
\end{lemma}

For a representation $A$ of
dimension $\vec n$ and each sequence
$C_1,\dots,C_t$ of $n_1\times n_1$,
$\dots$, $n_t\times n_t$ matrices,
we define the representation $[C,A]$
of the same dimension as follows:
\begin{equation}\label{eee}
[C,A]_{\alpha}=
C_qA_{\alpha}-A_{\alpha}C_p\qquad
\text{for each arrow $\alpha\colon
p\to q$}.
\end{equation}
Denote by $[{\mathbb F}^{\,\vec
n\times \vec n}, A]$ the set of such
representations.

Due to the next theorem, each
representation $A\in{\cal R}(\vec
n,\mathbb F)$ possesses a miniversal
deformation of the form \eqref{a2},
which was called in \cite{gar_ser} a
\emph{simplest miniversal
deformation} of $A$.

\begin{theorem}[{\cite[Theorem
2.1]{gar_ser}}]
 \label{t2.1}
Let $A$ be a matrix representation
of dimension $\vec n$ of a quiver
$Q$ with vertices $1,\dots,t$ over a
field $\mathbb{F}$ of complex or
real numbers. For each subset
$\varGamma$ of the set \eqref{a1},
the deformation ${\cal
U}_{\varGamma}(\vec {\varepsilon})$
defined in \eqref{a2} is miniversal
if and only if the vector space
${\cal R}(\vec n,\mathbb F)$ of all
representations of dimension $\vec
n$ decomposes into the direct sum
\begin{equation}\label{a4}
{\cal R}(\vec n,\mathbb F)=[{\mathbb
F}^{\,\vec n\times \vec n}, A]
\oplus {\cal E}_{\varGamma},
\end{equation}
in which ${\cal E}_{\varGamma}$
denotes the subspace spanned by all
elementary representations
$E_{\alpha ij}$ with
$(\alpha,i,j)\in\varGamma$.
\end{theorem}

In Section \ref{s3} we give a direct
proof of Theorem \ref{t2.1}. A
simplest miniversal deformation of
$A\in{\cal R}(\vec n,\mathbb F)$ can
be constructed as follows. Let
$T_1,\dots,T_r$ be a basis of the
space $[{\mathbb F}^{\,\vec n\times
\vec n}, A]$, and let
$E_1,\dots,E_l$ be the basis of
${\cal R}(\vec n,\mathbb F)$
consisting of all elementary
representations $E_{\alpha ij}$.
Removing from the sequence
$T_1,\dots, T_r,E_1,\dots,E_l$ every
representation that is a linear
combination of the preceding
representations, we obtain a new
basis $T_1,\dots, T_r,
E_{i_1},\dots,E_{i_k}$ of the space
$ {\cal R}(\vec n,\mathbb F)$. By
Theorem \ref{t2.1}, the deformation
$$ {\cal A}(\varepsilon_1,\dots,
\varepsilon_k)= A+\varepsilon_1
E_{i_1}+\dots+\varepsilon_kE_{i_k}
$$ is miniversal.

A \emph{direct sum} of two matrix
representations $A$ and $B$ is the
representation
\[
C=A\oplus B,\qquad
C_{\alpha}:=A_{\alpha}\oplus
B_{\alpha}\ \text{for all arrows
${\alpha}$.}
\]
A representation is called
\emph{indecomposable} if it is not
isomorphic to a direct sum of
representations of smaller sizes. It
is known that each matrix
representation $A$ is isomorphic to
a direct sum of indecomposable
representations
\begin{equation*}\label{a6}
A_1\oplus A_2\oplus\dots\oplus A_s
\end{equation*}
determined by $A$ uniquely up to
permutation of summands and
replacement them by isomorphic
representations.

\begin{theorem}\label{the}
Let $A=A_1\oplus\dots\oplus A_s$ be
a matrix representation, and let
\begin{multline}\label{a7}
{\cal A}(\vec
{\varepsilon})=A+{B}(\vec
{\varepsilon})\\=\begin{bmatrix}
  A_1+B_{11}(\vec
{\varepsilon}) & B_{12}(\vec
{\varepsilon}) &\dots & B_{1s}(\vec
{\varepsilon}) \\
  B_{21}(\vec
{\varepsilon}) &A_2+ B_{22}(\vec
{\varepsilon}) & \dots & B_{2s}(\vec
{\varepsilon}) \\
  \vdots & \vdots & \ddots & \vdots \\
  B_{s1}(\vec
{\varepsilon}) & B_{s2}(\vec
{\varepsilon}) & \dots &
A_s+B_{ss}(\vec {\varepsilon})
\end{bmatrix}
\end{multline}
be its deformation of the form
\eqref{a2}, whose matrices are
partitioned into blocks conformably
to the partition of $A$. Then ${\cal
A}(\vec {\varepsilon})$ is a
simplest mini\-versal deformation of
$A$ if and only if each
\begin{equation}\label{a8}
{\cal A}_{pq}(\vec
{\varepsilon}):=\begin{bmatrix}
  A_p+B_{pp}(\vec
{\varepsilon}) & B_{pq}\\
  B_{qp}(\vec
{\varepsilon}) &A_q+ B_{qq}(\vec
{\varepsilon})
\end{bmatrix},\qquad p<q,
\end{equation}
is a simplest miniversal deformation
of $A_p\oplus A_q$.
\end{theorem}

\begin{proof}
For $p,q\in\{1,\dots,s\}$ and for
each representation $M$ of dimension
$\vec n$, whose matrices are
partitioned into blocks conformably
to the partition of \eqref{a7},
denote by $M^{(p,q)}$ the
representation obtained from $M$ as
follows: in each of its matrices one
replaces by $0$ all blocks except
for the $(p,q)^{\text{th}}$ block.
If $\cal V$ is a  subspace of ${\cal
R}(\vec n,\mathbb F)$, then
\[
{\cal
V}^{(p,q)}:=\{M^{(p,q)}\,|\,M\in
{\cal V}\}
\]
is also a subspace.

Let the deformation \eqref{a7} be
miniversal. Since it has the form
\eqref{a2}, the decomposition
\eqref{a4} holds, and so each
$M\in{\cal R}(\vec n,\mathbb
F)^{(p,q)}$ is uniquely represented
in the form
\[
M=P+Q,\qquad P\in [{\mathbb
F}^{\,\vec n\times \vec n}, A],\quad
Q\in {\cal E}_{\varGamma}.
\]
Then
$M=M^{(p,q)}=P^{(p,q)}+Q^{(p,q)}$
and due to the obvious inclusions
\[
[{\mathbb F}^{\,\vec n\times \vec
n}, A]^{(p,q)}\subset [{\mathbb
F}^{\,\vec n\times \vec n},
A],\qquad {\cal
E}_{\varGamma}^{(p,q)}\subset {\cal
E}_{\varGamma}
\]
we have
\[
P\in [{\mathbb F}^{\,\vec n\times
\vec n}, A]^{(p,q)},\qquad Q\in
{\cal E}_{\varGamma}^{(p,q)},
\]
and so
\begin{equation}\label{a10}
{\cal R}(\vec n,\mathbb
F)^{(p,q)}=[{\mathbb F}^{\,\vec
n\times \vec n}, A]^{(p,q)} \oplus
{\cal E}_{\varGamma}^{(p,q)}
\end{equation}
for all $p,q\in\{1,\dots,s\}$. Due
to Theorem \ref{t2.1}, the
deformations \eqref{a8} are
miniversal for all $p<q$.

Conversely, if the deformations
\eqref{a8} are miniversal for all
$p<q$, then applying the same
reasoning as above to \eqref{a8}
instead of ${\cal A}(\vec
{\varepsilon})$, we obtain the
decompositions \eqref{a10} for all
$(p,q)$. They ensure the
decomposition \eqref{a4}, and so the
deformation \eqref{a7} is miniversal
by Theorem \ref{t2.1}.
\end{proof}

The next lemma helps to construct
miniversal deformations and will be
used in Section \ref{s2}.

\begin{lemma}\label{lems}
Let $A$ be a representation of $Q$
such that $A_{\alpha}=I$ for some
arrow $\alpha\colon p_1\to p_2$,
$p_1\ne p_2$. Denote by $Q'$ the
quiver obtained from $Q$ by removing
the arrow $\alpha$ and replacing
$p_1$ and $p_2$ by a single vertex
$p$ $($then each other arrow that
connects $p_1$ and $p_2$ becomes a
loop$)$. Denote by $A'$ the
representation of $Q'$ that is
obtained from $A$ by removing
$A_{\alpha}=I$. Then each miniversal
deformation of $A'$ can be extended
to a miniversal deformation of $A$
by assigning the identity matrix to
$\alpha$.
\end{lemma}

\begin{proof}
Let ${\cal
A}'(\lambda_1,\dots,\lambda_k)$ be a
miniversal deformation of $A'$, and
let ${\cal
A}(\lambda_1,\dots,\lambda_k)$ be
the deformation of $A$ obtained from
${\cal
A}'(\lambda_1,\dots,\lambda_k)$ by
assigning the identity matrix to
$\alpha$. We need to prove that
${\cal
A}(\lambda_1,\dots,\lambda_k)$
satisfies the condition (ii) of
Lemma \ref{lem}. Since
$A_{\alpha}=I$, the deformation
${\cal U}(\vec {\varepsilon})$ of
$A$ is equivalent to some
deformation ${\cal B}(\vec
{\varepsilon})$ that is the identity
on $\alpha$. Denote by ${\cal
B}'(\vec {\varepsilon})$ the
deformation of $A'$ obtained from
${\cal B}(\vec {\varepsilon})$ by
removing ${\cal B}(\vec
{\varepsilon})_{\alpha}=I$. Since
${\cal
A}'(\lambda_1,\dots,\lambda_k)$ is a
miniversal deformation of $A'$, by
Definition \ref{d} ${\cal B}'(\vec
{\varepsilon})$ is equivalent to a
deformation of the form ${\cal
A}'(\varphi_1(\vec\mu),\dots,
\varphi_k(\vec\mu)),$ where
$\varphi_i(\vec\mu)$ are convergent
in a neighborhood of $\vec 0$ power
series such that $\varphi_i(\vec
0)=0$. Then ${\cal B}(\vec
{\varepsilon})$ is equivalent to the
deformation ${\cal
A}(\varphi_1(\vec\mu),\dots,
\varphi_k(\vec\mu))$ and so ${\cal
A}(\lambda_1,\dots,\lambda_k)$
satisfies the condition (ii) of
Lemma \ref{lem}.
\end{proof}

\section{Miniversal deformations of
matrices of chains of linear
mappings}\label{s2}

In this section, we give simplest
miniversal deformations of matrices
of chains of linear mappings
$V_1-V_2-\dots-V_t$ over complex or
real numbers; that is, of
representations of the quiver
\begin{equation} \label{x1.5}
1\ \frac{{\alpha}_1}{\quad\qquad}\
  2\ \frac{{\alpha}_2}{\quad\qquad}\
  \cdots\
  \frac{{\alpha}_{t-2}}{\quad\qquad}\
  {(t-1)}
  \frac{{\alpha}_{t-1}}{\quad\qquad}\
  t
\end{equation}
in which each line denotes
$\longrightarrow$ or
$\longleftarrow$. Due to Theorem
\ref{the}, it suffices to give
simplest miniversal deformations of
those of its representations that
are direct sums of two
nonindecomposable representations.
Each representation $A$ of this
quiver is isomorphic to a direct
sum, determined uniquely up to
permutation of summands, of
indecomposable representations of
the form
\begin{equation}\label{x1.4}
L_{ij} :\qquad 1\ \frac{0}{\qquad}\
\cdots \ \frac{0}{\qquad}\ i\
\frac{I_{1}}{\qquad}\ \cdots\
\frac{I_{1}}{\qquad}\ j\
\frac{0}{\qquad}\ \cdots\
\frac{0}{\qquad}\ t\,,
\end{equation}
$1\le i\le j\le t$, having dimension
$(0,\dots,0,1,\dots,1,0,\dots 0)$
(in \cite{ser1} this direct sum is
constructed by $A$ using only
unitary transformations). Note that
the zero matrices in \eqref{x1.4}
have sizes $0\times 0$, $0\times 1$,
or $1\times 0$; it is agreed that
there exists exactly one matrix,
denoted by $0_{n0}$, of size
$n\times 0$ and there exists exactly
one matrix, denoted by $0_{0n}$, of
size $0\times n$ for every
nonnegative integer $n$; they
represent the linear mappings $0\to
{\mathbb F}^n$ and ${\mathbb F}^n\to
0$ and are considered as zero
matrices. Then $$ M_{pq}\oplus
0_{m0}=\begin{bmatrix}
  M_{pq} & 0 \\
  0 &0_{m0}
\end{bmatrix}=\begin{bmatrix}
  M_{pq}& 0_{p0} \\
  0_{mq}& 0_{m0}
\end{bmatrix}=\begin{bmatrix}
M_{pq} \\ 0_{mq}
\end{bmatrix}
$$ and $$ M_{pq}\oplus
0_{0n}=\begin{bmatrix}
  M_{pq} & 0 \\
  0 & 0_{0n}
\end{bmatrix}=\begin{bmatrix}
  M_{pq}& 0_{pn} \\
  0_{0q}& 0_{0n}
\end{bmatrix}=\begin{bmatrix}
   M_{pq} & 0_{pn}
\end{bmatrix}
$$ for every $p\times q$ matrix
$M_{pq}$.

The next theorem gives simplest
miniversal deformations of all
direct sums of two indecomposable
representations $L_{ij}$ of the
quiver \eqref{x1.5}. Using them and
Theorem \ref{the}, one can construct
a simplest miniversal deformation of
any representation decomposed into a
direct sum of indecomposable
representations.

\begin{theorem}\label{theo}
Let $L_{pq}$ and $L_{rs}$ $(p\le r)$
be two nonindecomposable
representations of the form
\eqref{x1.4} of the quiver
\eqref{x1.5} over complex or real
numbers. Then a miniversal
deformation of $A=L_{pq}\oplus
L_{rs}$ has at most $1$ parameter.
Moreover, it has no parameters
$($and hence coincides with $A)$ in
all the cases except for the next
cases, in which the representations
$A$ and their simplest miniversal
deformations ${\cal A}(\lambda)$ are
the following:

{\rm(i)} $L_{pq}\oplus L_{q+1,s}$
$(p\le q<s)$,
\begin{equation}\label{eq1}
\cdots\ \frac{}{\qquad}\
  q\ \frac{[\lambda]}{\quad\qquad}\
  q+1\
  \frac{}{\qquad}\ \cdots
\end{equation}

{\rm(ii)} $L_{pq}\oplus L_{qs}$
$(p<q<s)$,
\begin{gather}\label{eq2}
\cdots\ \frac{}{\qquad}\
  q-1\ \xrightarrow{\begin{bmatrix}
    1 \\0
  \end{bmatrix}}\
  q\  \xrightarrow{\begin{bmatrix}
    \lambda &1
  \end{bmatrix}}\ q+1\
  \frac{}{\qquad}\ \cdots
\\
\label{eq3} \cdots\ \frac{}{\qquad}\
  q-1\ \xleftarrow{\begin{bmatrix}
    1 &0
  \end{bmatrix}}\
  q\  \xleftarrow{\begin{bmatrix}
    \lambda \\1
  \end{bmatrix}}\ q+1\
  \frac{}{\qquad}\ \cdots
\end{gather}

{\rm(iii)} $L_{pq}\oplus L_{rs}$
$(p<r\le s<q)$,
\begin{gather*}\label{eq4}
\cdots\ \frac{}{\qquad}\
  r-1\ \xrightarrow{\begin{bmatrix}
    1 \\0
  \end{bmatrix}}\
  r\ \frac{I_2}{\qquad}\ \cdots \
  \frac{I_2}{\qquad}\ s\
   \xleftarrow{\begin{bmatrix}
   1\\ \lambda
  \end{bmatrix}}\ s+1\
  \frac{}{\qquad}\ \cdots
\\
\label{eq5} \cdots\ \frac{}{\qquad}\
  r-1\ \xleftarrow{\begin{bmatrix}
    1 &0
  \end{bmatrix}}\
  r\ \frac{I_2}{\qquad}\ \cdots \
  \frac{I_2}{\qquad}\ s\
   \xrightarrow{\begin{bmatrix}
   1& \lambda
  \end{bmatrix}}\ s+1\
  \frac{}{\qquad}\ \cdots
\end{gather*}
$($each line denotes
$\longrightarrow$ or
$\longleftarrow$; all unspecified
matrices of ${\cal A}(\lambda)$
coincide with the corresponding
matrices of $A)$.
\end{theorem}

\begin{proof}
Let us find a miniversal deformation
of $A=L_{pq}\oplus L_{rs}$. We may
suppose that the pairs $(p,q)$ and
$(r,s)$ are lexicographically
ordered; that is, $p\le r$ and if
$p=r$ then $q\le s$. Deleting arrows
of the quiver \eqref{x1.5} that
correspond to matrices without rows
or columns, we reduce our
consideration to the case $p=1$,
$\max(q,s)=t$, and $r\le q+1$. Due
to Lemma \ref{lems}, we may suppose
that $A$ has no identity matrices.
Then the quiver has at most 3
vertices and $A$ is one of the
following representations:
\begin{equation*}\label{q1}
L_{11}\oplus L_{11},\ \ L_{11}\oplus
L_{22},\ \ L_{12}\oplus L_{22},\ \
L_{11}\oplus L_{12},\ \ L_{12}\oplus
L_{23},\ \ L_{13}\oplus L_{22},
\end{equation*}
or diagrammatically
\begin{equation*}\label{q2}
\begin{bmatrix}
   \bullet\\ \bullet
  \end{bmatrix},
       \ \
  \begin{bmatrix}
   \bullet &\\ & \bullet
  \end{bmatrix},
      \ \
  \begin{bmatrix}
   \bullet \frac{}{\quad} \bullet\\
   \phantom{\bullet \frac{}{\quad}} \bullet
  \end{bmatrix},
      \ \
 \begin{bmatrix}
   \bullet \phantom{\frac{}{\quad} \bullet}\\
   {\bullet \frac{}{\quad}} \bullet
  \end{bmatrix},
      \ \
\begin{bmatrix}
   \bullet \frac{}{\quad} \!\bullet\phantom{\! \frac{}{\quad} \bullet}\\
   \phantom{\bullet \frac{}{\quad} \!}\bullet\! \frac{}{\quad} \bullet
  \end{bmatrix},\ \
  \begin{bmatrix}
   \bullet \frac{}{\quad} \!\bullet{\! \frac{}{\quad} \bullet}\\
\bullet
  \end{bmatrix};
\end{equation*}
here the first row of each matrix
represents $L_{pq}$ and the second
row represents $L_{rs}$.

The representation $L_{11}\oplus
L_{11}$ has no matrices, and so its
deformation has no parameters.

The representation $L_{11}\oplus
L_{22}$ consists of the $1$-by-$1$
matrix $[0]$ and has the deformation
\eqref{eq1}.

The representation $L_{12}\oplus
L_{22}$ is
\[
  1\ \xrightarrow{\begin{bmatrix}
    1 \\0
  \end{bmatrix}}\
  2\qquad\text{or}\qquad
  1\  \xleftarrow{\begin{bmatrix}
    1 &0
  \end{bmatrix}}\ 2
\]
depending on the orientation of the
arrow. In both the cases, $[{\mathbb
F}^{\,\vec n\times \vec n}, A]={\cal
R}(\vec n,\mathbb F)$ and so by
\eqref{a4} ${\cal E}_{\varGamma}=0$.
Hence each miniversal deformation of
$A$ has no parameters. The same
holds for $A=L_{11}\oplus L_{12}$.

Depending on the orientation of
arrows, $A=L_{12}\oplus L_{23}$ has
one of the forms:
\begin{gather}\label{eqv5}
  1\ \xrightarrow{\begin{bmatrix}
    1 \\0
  \end{bmatrix}}\
  2\  \xrightarrow{\begin{bmatrix}
    0 &1
  \end{bmatrix}}\ 3
     \qquad\qquad
  1\ \xleftarrow{\begin{bmatrix}
    1 &0
  \end{bmatrix}}\
  2\  \xleftarrow{\begin{bmatrix}
    0 \\1
  \end{bmatrix}}\ 3
             \\ \label{eqv8}
  1\ \xrightarrow{\begin{bmatrix}
    1 \\0
  \end{bmatrix}}\
  2\  \xleftarrow{\begin{bmatrix}
    0 \\1
  \end{bmatrix}}\ 3
     \qquad\qquad
  1\ \xleftarrow{\begin{bmatrix}
    1 &0
  \end{bmatrix}}\
  2\  \xrightarrow{\begin{bmatrix}
    0 &1
  \end{bmatrix}}\ 3
  \end{gather}
In the first case, the space
$[{\mathbb F}^{\,\vec n\times \vec
n}, A]$ with $\vec n=(1,2,1)$
consists of the representations
\[
[C,A]\colon\qquad 1\
\xrightarrow{\begin{bmatrix}
    b_1-a \\b_3
  \end{bmatrix}}\
  2\  \xrightarrow{\begin{bmatrix}
    -b_3 &c-b_4
  \end{bmatrix}}\ 3
\]
defined by \eqref{eee}, in which
\[
C_1=[a],\qquad C_2=\begin{bmatrix}
    b_1 &b_2\\    b_3 &b_4
  \end{bmatrix},\qquad C_3=[c].
\]
Due to \eqref{a4}, the representation
\eqref{eq2} is a miniversal deformation
of $A$. Analogously, \eqref{eq3} is a
miniversal deformation of the second
representation in \eqref{eqv5}. If $A$
is of the form \eqref{eqv8}, then
$[{\mathbb F}^{\,\vec n\times \vec
n},A]$ coincides with ${\cal R}(\vec
n,\mathbb F)$ and so each miniversal
deformation of $A$ has no parameters.

Depending on the orientation of
arrows, $A=L_{13}\oplus L_{22}$ has
one of the forms:
\begin{gather}\label{eqv51}
  1\ \xrightarrow{\begin{bmatrix}
    1 \\0
  \end{bmatrix}}\
  2\  \xrightarrow{\begin{bmatrix}
    1&0
  \end{bmatrix}}\ 3
     \qquad\qquad
  1\ \xleftarrow{\begin{bmatrix}
    1 &0
  \end{bmatrix}}\
  2\  \xleftarrow{\begin{bmatrix}
    1 \\0
  \end{bmatrix}}\ 3
             \\ \label{eqv81}
  1\ \xrightarrow{\begin{bmatrix}
    1 \\0
  \end{bmatrix}}\
  2\  \xleftarrow{\begin{bmatrix}
    1 \\0
  \end{bmatrix}}\ 3
     \qquad\qquad
  1\ \xleftarrow{\begin{bmatrix}
    1 &0
  \end{bmatrix}}\
  2\  \xrightarrow{\begin{bmatrix}
    1 &0
  \end{bmatrix}}\ 3
  \end{gather}
All miniversal deformations of
\eqref{eqv51} have no parameters.
The representations (iii) in Theorem
\ref{theo} are miniversal
deformations of \eqref{eqv81}.
\end{proof}

\section{A direct proof of Theorem
\ref{t2.1}} \label{s3}

For each matrix $P=[p_{ij}]$ over a
field $\mathbb F$ of complex or real
numbers, we define its \emph{norm}
as follows:
\begin{equation*}\label{4a}
\|P\|:=\sum |p_{ij}|.
\end{equation*}
By \cite[Section 5.6]{hor_John},
\begin{equation}\label{lk}
\|aP+bQ\|\le
|a|\,\|P\|+|b|\,\|Q\|,\qquad
\|PQ\|\le \|P\|\,\|Q\|
\end{equation}
for matrices $P$ and $Q$ and
$a,b\in\mathbb F$.

For each finite set
${M}=\{M_1,\dots,M_l\}$ of matrices,
we put
\[
\|M\|:=\|M_1\|+\dots+\|M_l\|.
\]

Let $Q$ be a quiver with vertices
$1,\dots,t$, let $M$ be its
representation of dimension $\vec
n=(n_1,\dots,n_t)$, and let
$S=(S_1,\dots,S_t)$ be a sequence of
matrices of sizes $n_1\times
n_1,\dots,n_t\times n_t$ (such
sequences will be called $\vec
n$-\emph{sequences}; they are closed
under addition and multiplication).
Denote by $SM$ and $MS$ the
representations of $Q$ obtained from
$M$ by replacing each matrix
$M_{\alpha}$ assigned to
$\alpha\colon p\to q$ with $S_q
M_{\alpha}$ and, respectively,
$M_{\alpha}S_p$. Due to \eqref{lk},
\[
\|SM\|\le\sum_{p,\alpha}
\|S_p\|\,\|M_{\alpha}\|=\|S\|\,\|M\|,
\qquad \|M\|\,\|S\|\le\|M\|\,\|S\|.
\]
If $M_{\alpha}=[m_{\alpha ij}]$ and
$\varGamma\subset\Upsilon_{\vec n}$
is a subset of \eqref{a1}, then we
put
\begin{equation*}\label{5aaa}
\|M\|_{\varGamma}:=
\sum_{(\alpha,i,j)\notin{\varGamma}}
|m_{\alpha ij}|;
\end{equation*}
in particular,
$\|M\|_{\Upsilon_{\vec n}}=\|M\|$.

\begin{lemma} \label{lem1az}
Let $A$ and $\varGamma$ be the
representation and the set from
Theorem {\rm\ref{s3}}. There exists
a natural number $m$ such that for
each real numbers $\varepsilon$ and
$\delta$ satisfying
\begin{equation*}\label{16z}
0<\varepsilon\le \delta<\frac 1m
\end{equation*}
and for each representation $M$ of
$Q$ satisfying
\begin{equation}\label{15z}
\|M\|_{\varGamma}<\varepsilon,\quad
\|M\|<\delta
\end{equation}
there exists an $\vec n$-sequence
\begin{equation}\label{5az}
S=I_{\vec n}+X,
\qquad\|X\|<m\varepsilon,
\end{equation}
in which $I_{\vec n}
=(I_{n_1},\dots,I_{n_t})$ and the
entries of matrices of $X$ are
linear polynomials in entries of $M$
such that
\begin{equation}\label{7z}
S(A+M)S^{-1}=A+M',\qquad
\|M'\|_{\varGamma}
<m\varepsilon\delta,\quad
\|M'\|<\delta+m\varepsilon.
\end{equation}
\end{lemma}

\begin{proof}
First we construct the $\vec
n$-sequence \eqref{5az}. By
\eqref{a4}, for each elementary
representation $E_{\alpha ij}$,
$(\alpha,i,j)\in\Upsilon_{\vec n} $,
(they were introduced after
Definition \ref{d}), there exists an
$\vec n$-sequence $X_{\alpha ij}$
such that
\begin{equation*}\label{8z}
E_{\alpha ij}+X_{\alpha
ij}A-AX_{\alpha ij}\in{\cal
E}_{\varGamma}.
\end{equation*}
If $M=\sum_{\alpha ij} m_{\alpha
ij}E_{\alpha ij}$ (that is, the
representation $M$ from Lemma
\ref{lem1az} is formed by the
matrices $M_{\alpha} =[m_{\alpha
ij}]$), then
\[
\sum_{\alpha ij} m_{\alpha
ij}E_{\alpha ij}+\sum_{\alpha ij}
m_{\alpha ij} X_{\alpha
ij}A-\sum_{\alpha ij} m_{\alpha
ij}AX_{\alpha ij}\in {\cal
E}_{\varGamma}
\]
and for
\[
S=I_{\vec n}+X,\qquad
X:=\sum_{\alpha ij} m_{\alpha
ij}X_{\alpha ij},
\]
we have
\begin{equation}\label{18z}
M+SA-AS\in {\cal E}_{\varGamma}.
\end{equation}
If $(\alpha, i,j)\in{\varGamma}$,
then $E_{\alpha ij}\in {\cal
E}_{\varGamma}$ and we can put
$X_{\alpha ij}=0$. If $(\alpha,
i,j)\notin{\varGamma}$, then
$|m_{\alpha ij}|<\varepsilon$ by the
first inequality in \eqref{15z}. We
obtain
\begin{equation}\label{jio}
\|X\|\le \sum_{(\alpha,
i,j)\notin{\varGamma}} |m_{\alpha
ij}|\,\|X_{\alpha ij}\|<
\sum_{(\alpha,
i,j)\notin{\varGamma}}
\varepsilon\|X_{\alpha ij}\|=
\varepsilon c,
\end{equation}
where
\[
\qquad c:= \sum_{(\alpha,
i,j)\notin{\varGamma}} \|X_{\alpha
ij}\|.
\]

Take $\varepsilon<1/(2c)$, then
\begin{equation}\label{12z}
\varepsilon c<\frac 1 2
\end{equation}
and so
\[
\|X^k\|\le \|X\|^k<(\varepsilon
c)^k<1/2^k\to 0\qquad\text{if $k\to
\infty$.}
\]
Hence,
\begin{align}\nonumber
S^{-1}&=(I_{\vec n}+X)^{-1}=
I_{\vec n}-X+X^2-X^3+\cdots\\
\label{10z} &=I_{\vec
n}-XS^{-1}=I_{\vec n}-X+X^2S^{-1}.
\end{align}
Furthermore,
\[
\|S^{-1}\|\le \|I_{\vec
n}\|+\|X\|+\|X\|^2+\cdots<n+
\varepsilon c+(\varepsilon
c)^2+\cdots,
\]
where $n:=n_1+\dots+n_t$, and by
\eqref{12z}
\begin{equation}\label{13z}
\|S^{-1}\|\le n-1+\frac
1{1-1/2}=n+1.
\end{equation}

Using \eqref{10z}, we obtain
\begin{multline*}
S(A+M)S^{-1} =(A+M+XA+XM)S^{-1}
=A(I_{\vec n}-X+X^2S^{-1})\\
+(M+XA)(I_{\vec n}-XS^{-1})+XMS^{-1}
=A+M',
\end{multline*}
where
\begin{align*}\label{14z}
M'&:=M+XA-AX+N,\\
N&:=AX^2S^{-1}-(M+XA)XS^{-1}+XMS^{-1}.
\end{align*}
Then by \eqref{jio}, \eqref{13z},
and \eqref{15z}, and since
$\varepsilon\le \delta$, we have
\begin{align*}
\|N\|\le 2\|A\|(\varepsilon
c)^2(n+1)+2\delta(\varepsilon
c)(n+1)\le \varepsilon\delta d,
\end{align*}
where
\[
d:=2\|A\| c^2(n+1)+2c(n+1).
\]

By \eqref{18z}, $M+XA-AX\in {\cal
E}_{\varGamma}$, and so
\[
\|M'\|_{\varGamma}=\|N\|_{\varGamma}\le
\varepsilon\delta d.
\]
Furthermore,
\[
\|M'\|\le \|M\|+\|XA-AX\|+\|N\|\le
\delta+2\varepsilon c \|A\|
+\varepsilon\delta d= \delta +e
\varepsilon,
\]
where
\[
e:=2c\|A\|+\delta d.
\]
Taking any natural number $m$ that
is greater than $c$, $d$, and $e$,
we obtain \eqref{5az} and
\eqref{7z}.
\end{proof}

\begin{lemma} \label{lem2z}
Let $m$ be any natural number being
$\ge 3$, and let
\begin{equation*}\label{21z}
\varepsilon_1,\ \delta_1,\
\varepsilon_2,\ \delta_2,\
\varepsilon_3,\ \delta_3,\,\dots
\end{equation*}
be the sequence of numbers defined
by induction:
\begin{equation}\label{22z}
\varepsilon_1=\delta_1=m^{-7},\quad
\varepsilon_{i+1}=m\varepsilon_i
\delta_i,\quad
\delta_{i+1}=\delta_i+m\varepsilon_i.
\end{equation}
Then
\begin{equation}\label{23z}
\varepsilon_{i}<m^{-4i},\quad
\delta_i<m^{-5}
\end{equation}
for all $i$ and
\begin{equation}\label{24z}
\varepsilon_1+\varepsilon_2+
\varepsilon_3+\dots<2.
\end{equation}
\end{lemma}

\begin{proof}
Reasoning by induction, we assume
that the inequalities \eqref{23z}
hold for $i=1,2,\dots,l$. Then
\[
\varepsilon_{l+1}=m\varepsilon_l
\delta_l<m^{-4l}mm^{-5}=m^{-4(l+1)}
\]
and
\begin{align*}
\delta_{l+1}&=\delta_l+m\varepsilon_l=
\delta_{l-1}+m(\varepsilon_{l-1}+
\varepsilon_l)=\cdots\\&=
\delta_{1}+m(\varepsilon_{1}
+\varepsilon_{2}+\dots+
\varepsilon_l)\\
&<m^{-7}+m(m^{-7}+m^{-4\cdot 2}
+m^{-4\cdot 3}+\cdots)\\
 &\le
m^{-7}+m^{-6}(1+m^{-1}+m^{-2}
+m^{-3}+\cdots)
 \\&
=m^{-7}+m^{-6}\frac 1{1-m^{-1}}\le
m^{-6}\left(m^{-1}+\frac
32\right)<2m^{-6}<m^{-5}.
\end{align*}
This proves \eqref{23z} for all $i$.
Then \eqref{24z} holds too since
\[
\varepsilon_1+\varepsilon_2+
\varepsilon_3+\dots<
m^{-4}+m^{-4\cdot 2} +m^{-4\cdot
3}+\cdots<\frac 1{1-m^{-4}}<2.
\]
\end{proof}

Theorem \ref{t2.1} follows from the
next lemma.

\begin{lemma}\label{lem3z}
Let $A$ and $\varGamma$ be from Theorem
{\rm\ref{s3}}, let $m$ be a natural
number that is greater than $3$ and
satisfies Lemma \ref{lem1az}, and let
$M$ be any representation satisfying
$\|M\|<m^{-7}$. Then there exists an
$\vec n$-sequence $S=I_{\vec n}+X$
depending holomorphically on the
entries of $M$ in a neighborhood of
zero such that $$S(A+M)S^{-1}-A\in
{\cal E}_{\varGamma}$$ and $S=I_{\vec
n}$ if $M=0$.
\end{lemma}

\begin{proof}
We construct a sequence of
representations
\begin{equation*}\label{20aaz}
A+M_1,\ A+M_2,\ A+M_3,\dots
\end{equation*}
by induction. Put $M_1=M$. Let $M_i$
be constructed and let
\begin{equation*}\label{15nz}
\|M_i\|_{\varGamma}<\varepsilon_i,\quad
\|M_i\|<\delta_i,
\end{equation*}
where $\varepsilon_i$ and $\delta_i$
are defined in \eqref{22z}. Then by
\eqref{23z} and Lemma \ref{lem1az}
there exists
\begin{equation}\label{5aaz}
S_{i+1}=I_{\vec n}+X_{i+1},\qquad
\|X_{i+1}\|<m\varepsilon_{i+1},
\end{equation}
such that
\begin{equation*}\label{7aaz}
S_{i+1}(A+M_{i})S_{i+1}^{-1}=A+M_{i+1},
\quad \|M_{i+1}\|_{\varGamma}
<\varepsilon_{i+1},\quad
\|M_{i+1}\|<\delta_{i+1}.
\end{equation*}
For each natural number $l$, put
\begin{equation}\label{28z}
S^{(l)}:=S_l\cdots S_3S_2S_1 =(I_{\vec
n}+X_l)\cdots(I_{\vec n}+X_2)(I_{\vec
n}+X_1).
\end{equation}
Let us prove that the sequence
$S^{(1)},S^{(2)},S^{(3)},\dots$
converges. Indeed,
\[
S^{(l)}=I_{\vec n}+\sum_{l\ge
i}X_i+\sum_{l\ge i>j}X_iX_j+\cdots
\]
and so
\begin{align}\nonumber
\|S^{(l)}\|&\le\|I_{\vec
n}\|+\sum_{l\ge i}\|X_i\|
+\sum_{l\ge i>j}\|X_i\|\,\|X_j\|+\cdots\\
\label{lo} &\le
n-1+(1+\|X_1\|)(1+\|X_2\|)
(1+\|X_3\|)\cdots,
\end{align}
where $n:=n_1+\dots+n_t$. By
\cite[Section III, \S 4.3]{mark},
the product \eqref{lo} converges
since the sum
\[
\|X_1\|+\|X_2\|+\|X_3\|+\cdots
\]
converges due to \eqref{5aaz} and
\eqref{24z}.

The entries of all matrices forming
$S:=\lim S^{(l)}$ are holomorphic
functions in the entries of $M$
(that satisfies $\|M\|<m^{-7}$) due
to Weierstrass' theorem
\cite[Section III, \S 4.1]{mark}
since the sequence
\begin{equation*}\label{29z}
I_{\vec n}+X_1,\ (I_{\vec
n}+X_2)(I_{\vec n}+X_1),\ (I_{\vec
n}+X_3)(I_{\vec n}+X_2)(I_{\vec
n}+X_1),\,\dots
\end{equation*}
converges uniformly to \eqref{28z}.

Since $A+M_l\to S(A+M)S^{-1}$ if
$l\to\infty$ and
$\|M_{l}\|_{\varGamma}
<\varepsilon_{l}\to 0$, we have
$S(A+M)S^{-1}-A\in {\cal
E}_{\varGamma}$.
\end{proof}

\end{document}